\documentclass[twoside,12pt]{amsart} 
\usepackage[dvips]{graphicx}
\DeclareGraphicsExtensions{ps,eps}
\usepackage{amsfonts,mathrsfs}
\theoremstyle{plain}
\begingroup
\newtheorem{theorem}{Theorem}[section]

\endgroup

\theoremstyle{definition}
\begingroup

\endgroup

\theoremstyle{remark}
\begingroup
\endgroup

\mathchardef\emptyset="001F

\numberwithin{equation}{section}

\newcommand{\ERRE}{{{\mathbb{R}}}}
\newcommand{\ZED}{{{\mathbb{Z}}}}
\newcommand{\LERAY}{{{\mathbb{P}}}}

\title[Combinatorics, Optimal Transport and Hydrodynamics
]{
Connections between Optimal Transport, Combinatorial Optimization and Hydrodynamics
}

\author[Yann Brenier, CNRS, CMLS]{Yann Brenier, CNRS, 
Centre de Math\'ematiques Laurent Schwartz, 
Ecole Polytechnique, Palaiseau, France}

\begin{document}
\begin{abstract}

There are well-established
connections between combinatorial optimization, optimal transport theory
and Hydrodynamics, through the linear assignment problem in combinatorics, the
Monge-Kantorovich problem in optimal transport theory and the model of inviscid,
pressure-less fluids in Hydrodynamics.
Here, we consider the more challenging quadratic assignment problem (which is NP,
while the linear assignment problem is just P)
and find, in some particular case,
a correspondence with the problem of finding 
stationary solutions of Euler's equations for incompressible
fluids. For that purpose, we introduce and analyze a "gradient flow" equation
\begin{equation}
\label{}
\;\;\;\
\partial_t \varphi+\nabla\cdot(\varphi v)=0,\;\;\;(-\bigtriangleup)^m v=-\LERAY\nabla\cdot(\nabla\varphi\otimes\nabla\varphi),
\end{equation}
where ${\LERAY}$ denotes the ${L^2}$ projection onto divergence-free vector fields
and $m=0$ or $m=1$, with suitable boundary conditions.
Then,
combining some ideas of P.-L. Lions (for the Euler equations) and Ambrosio-Gigli-Savar\'e (for the heat equation),
we provide for the initial value problem a 
concept
of generalized ``dissipative'' solutions which always exist globally in time and are unique whenever they
are smooth.

\end{abstract}
\maketitle

\section{Well-known connections between Optimal transport theory, Hydrodynamics
and combinatorial optimization}
\subsection{The Monge-Kantorovich distance in optimal transport theory}
The (quadratic) Monge-Kantorovich ($MK^2$) distance 
(very often called ``Wasserstein'' distance in optimal transport theory and also called Tanaka distance in kinetic theory \cite{Vi}) can be defined in terms of probability measures and random variables as:
\begin{equation}
\label{}
d_{MK^2}(\mu,\nu)
=\;\;\inf\{\sqrt{E(|X-Y|^2)},\;\;\;{\rm{
{law}}}(X)=\mu,\;{\rm{
{law}}}(Y)=\nu\}
\end{equation}
where $\mu$ and $\nu$ are probability measures (with finite second moments) defined on the Euclidean space ${\mathbf{R}^d}$, $X$ and $Y$ denotes random variables valued in ${\mathbf{R}^d}$, $|\cdot|$  is the Euclidean norm and
$E$ denotes the expected value.

\subsection{Hydrodynamic interpretation of the $MK^2$ distance}

Using the so-called "Benamou-Brenier formula" or the "Otto calculus" \cite{BeBr,Ot,AGS1}, 
we may express the
$MK^2$ distance in hydrodynamic terms. More precisely, at least in the case when $\mu$ and $\nu$
are absolutely continuous with respect to the Lebesgue measure,
we may write
\begin{equation}
\label{}
d_{MK^2}^2(\mu,\nu)=
\inf\;\;\int_0^1\int_{\mathbf{R}^d}|v(t,x)|^2\rho(t,x)dtdx
\end{equation}
where the infimum is taken over all density and velocity fields 
$$
(t,x)\in [0,1]\times \mathbf{R}^d\rightarrow (\rho(t,x),v(t,x))\in \mathbf{R}_+ \times \mathbf{R}^d
$$
subject to the continuity equation
\begin{equation}
\label{}
\partial_t\rho+\nabla_x\cdot (\rho v)=0
\end{equation}
and the time-boundary conditions
\begin{equation}
\label{}
\rho(t=0,x)dx=\mu(dx),\;\;\rho(t=1,x)dx=\nu(dx).
\end{equation}
The (formal) optimality equations read
\begin{equation}
\label{}
v(t,x)=\nabla_x\theta(t,x),\;\;\; 
\partial_t\theta+\frac{1}{2}|\nabla_x\theta|^2=0,
\end{equation}
and describe a potential, inviscid, pressure-less gas, 
sometimes called "dust" (in cosmology in particular), which is one of the most trivial models of fluids.

\subsection{$MK^2$ distance and combinatorial optimization
}
Given two discrete probability measures 
on ${\mathbf{R}^d}$
\begin{equation}
\label{}
\mu=\sum_{i=1,N}\delta _{A_i},\;\;\;\nu=\sum_{j=1,N}\delta _{B_{j}},
\end{equation}
we easily check that
\begin{equation}
\label{}
d_{MK^2}^2(\mu,\nu)
=\inf_{{\rm{{law}}}(X)=\mu,{\rm{{law}}}(Y)=\nu}E(|X-Y|^2)
=\inf_{\sigma\in\mathcal{S}_N}
\sum_{i=1,N}\frac{|A_i-B_{\sigma_i}|^2}{N},
\end{equation}
where $\sigma\in\mathcal{S}_N$ denotes the set of all permutations of
the first $N$ integers.
Thus, computing the $MK^2$ distance between two discrete measures
is equivalent to solving the so-called
"linear assignment problem" (LAP):
\begin{equation}
\label{}
\inf_{\sigma\in\mathcal{S}_N}\; \sum_{i=1}^N c(i,\sigma_i),
\end{equation}
in the special case when the "cost matrix" $c$ has geometric contain
\begin{equation}
\label{}
c(i,j)=|A_{i}-B_{j}|^2.
\end{equation}
In full generality, the LAP is one of the simplest combinatorial optimization problems,
with complexity ${O(N^3)}$ \cite{Ba}.

\section{NP combinatorial optimization problems and Hydrodynamics}
There are much more challenging problems in combinatorial optimization,
such as the (NP) "quadratic assignment problem"
(which includes the famous traveling salesman problem).
\\
Given two ${N\times N}$ matrices ${\gamma}$ and ${c}$, 
with coefficients $\ge 0$, solve:
\begin{equation}
\label{QAP}
\;\;(QAP)\;\;\;\inf_{\;\sigma\in\mathcal{S}_N}\sum_{i,j=1,N}c(\sigma_i,\sigma_j) \gamma(i,j).
\end{equation}
The QAP is useful in computer vision \cite{Me}. Some continuous versions of the QAP are 
related to recent works
in geometric and functional analysis \cite{St}.

It turns out that the QAP can also be related to Hydrodynamics, as we are going to see.

\subsection{A minimization problem in Hydrodynamics}
This problem goes back to Lord Kelvin and has been frequently studied since
(see, for instance \cite{Be,Bu}).
Let
${D}$ be a smooth domain of unit Lebesgue measure in ${\mathbf{R}^d}$ and a 
real function $\varphi_0$ belonging to the space  $H^1_0(D)$ of Sobolev functions
vanishing along $\partial D$. We denote by
$\lambda$ the law of $\varphi_0$ over ${\mathbf{R}}$, so that
$$
\int_{D} F(\varphi_0(x))dx=\int_{-\infty}^\infty F(r)\lambda(dr),
$$
for all bounded continuous function
${F: \mathbf{R}\rightarrow \mathbf{R}}$.
We want to minimize the Dirichlet integral 
\begin{equation}
\label{dirichlet}
\mathcal{E}[\varphi]=\frac{1}{2}\int_{D} |\nabla\varphi(x)|^2dx
\end{equation}
among all real valued functions  
$\varphi\in H^1_0(D)$
with law $\lambda$, which may be written:
\begin{equation}
\label{kelvin}
\;\;\;\
\inf\{\; \frac{1}{2}\int_{D} |\nabla\varphi(x)|^2dx,\;\;
\varphi\in H^1_0(D),\;\;\;{\rm{Law}}(\varphi)={\rm{Law}}(\varphi_0)=\lambda\;\}
\end{equation}
or rephrased as a saddle-point problem:
\begin{equation}
\label{kelvin saddle}
\;\;\;\
\inf_{\varphi\in H^1_0(D)}\;\;\sup_{F: \mathbf{R}\rightarrow \mathbf{R}}
\;\;\;\;\; \frac{1}{2}\int_{D} |\nabla\varphi(x)|^2dx
+\int_{D} F(\varphi(x))dx-\int_{-\infty}^\infty F(r)\lambda(dr).
\end{equation}
Optimal solutions are formally solutions to
\begin{equation}
\label{}
\;\;
-\bigtriangleup \varphi+F'(\varphi)=0,
\;\;\;\varphi\in H^1_0(D),
\end{equation}
for some function ${F: \mathbf{R}\rightarrow \mathbf{R}}$, and, 
in 2d, are just stationary solutions to the Euler equations of incompressible fluids \cite{ArKh,MaPu}. (More precisely, $\varphi$ is the stream-function of
a stationary two-dimensional incompressible inviscid fluid.)
\subsection{The discrete version of the hydrodynamic problem is a QAP}

Let us discretize the domain $D$ with a lattice of $N$ vertices
$A_1,\cdot\cdot\cdot,A_N$ and define coefficients $\gamma(i,j)\ge 0$
so that the Dirichlet integral of a function $\varphi$ can be approximated
as follows:
$$
\int_D|\nabla\varphi(x)|^2dx\;\sim \;
\sum_{i,j=1}^N \gamma(i,j)|\varphi(A_i)-\varphi(A_j)|^2.
$$
At the discrete level, we may say that $\varphi$ and $\varphi_0$ have the
same (discrete law) whenever
$$
\varphi(A_i)=\varphi_0(A_{\sigma_i}),\;\;\;i=1,\cdot\cdot\cdot,N,
$$
for some permutation $\sigma\in\mathcal{S}_N$.
Thus, the discrete version of (\ref{kelvin}) reads:
\\
Find a
permutation ${\sigma}$ that achieves
\begin{equation}
\label{}
\;\;\;\
\inf_\sigma\;\;\; \sum_{i,j=1}^N c(\sigma_i,\sigma_j) \gamma(i,j)
\end{equation}
with
\begin{equation}
\label{}
\;\;\;\
c(i,j)=|\varphi_0(A_{i})-\varphi_0(A_{j})|^2.
\end{equation}
So we have clearly obtained a particular case of QAP
(\ref{QAP}).

\section{A "gradient-flow" approach to the hydrodynamic problem}

To address problem (\ref{kelvin}), it is natural to use a "gradient flow" approach involving a time dependent
function $\varphi_t(x)$ starting from $\varphi_0(x)$ at $t=0$. Hopefully, as $t\rightarrow +\infty$,
$\varphi_t$ will reach a solution to our minimization problem.
A canonical way of preserving the law $\lambda$ of
$\varphi(t,\cdot)$ during the evolution
is the transport of $\varphi$ by a (sufficiently) smooth
time-dependent divergence-free velocity field $v=v_t(x)\in \mathbf{R}^d$, parallel to $\partial D$,
according to
\begin{equation}
\label{transport}
\partial_t \varphi_t+\nabla\cdot(v_t\varphi_t)=0,\;\;\nabla\cdot v_t=0,\;\;\;v_t//\partial D.
\end{equation}
Indeed, we easily get:
$$
\frac{d}{dt}\int_D F(\varphi_t(x))dx
=-\int_D F'(\varphi_t(x))\nabla\cdot(v_t(x)\varphi_t(x))dx
$$
$$
=-\int_D v_t(x)\cdot\nabla(F(\varphi_t(x)))dx=0
$$
(since $v$ is divergence-free), for all smooth bounded function $F$.
Loosely speaking, the vector field $v$ should be interpreted as a 
kind of ``tangent vector''
along the ``orbit'' of all $\varphi$ sharing the same law $\lambda$
as $\varphi_0$.
\\
From the analysis viewpoint, according to the DiPerna-Lions theory on ODEs \cite{DiLi},
for the law $\lambda$ to be preserved, there is no need for $v$ to be very smooth
and it is just enough
that the space derivatives of $v$ 
are Lebesgue integrable functions (or even bounded Borel measures, according 
to Ambrosio \cite{Am}):
$$
\int_0^T\int_D |\nabla v_t(x)|dxdt<+\infty,\;\;\;\forall T>0.
$$
(N.B. in that situation, the solution of (\ref{transport})
is just (implicitly) given by
$
\varphi_t(\xi_t(x))=\varphi_0(x),
$
where $\xi$ is the unique time-dependent family of 
almost-everwhere one-to-one volume-preserving Borel maps
of $D$ generated by $v$ through:
$$
\partial_t \xi_t(x)=v_t(\xi_t(x)),\;\xi_0(x)=x.
$$
Of course, these maps are orientation preserving diffeomorphisms whenever
$v$ is smooth.)
\\
In order to get a ``gradient flow'', 
we also need a quadratic form (or, more generally, a convex functional,
which would then rather correspond to a ``Finslerian flow'') acting on
the ``tangent vector'' $v$.
For this purpose,
let us first denote by $\rm{Sol}(D)$ the Hilbert space of all square-integrable 
divergence free vector fields on $D$ and parallel to $\partial D$, with
$L^2$ norm and inner-product respectively denoted by $||\cdot||$
and $((\cdot,\cdot))$.
Next, let us fix a lower semi-continuous convex functional 
$K: a\in {\rm{Sol}(D)}\rightarrow K[a]\in [0,+\infty]$. We assume that 
at each smooth vector-field $\omega$ in $\rm{Sol}(D)$, 
\\
i) $K$ is finite,
\\
ii) its subgradient has a unique element $K'[\omega]\in L^2(D,\mathbf{R}^d)$, 
\\
iii) there is $\epsilon_K[\omega]>0$ such that, the "relative entropy" of $K$ controls the $L^2$ distance:
\begin{equation}
\label{entropy}
\eta_K[v,\omega]=K[v]-K[\omega]-((K'[\omega],v-\omega))\ge \epsilon_k[\omega]||v-\omega||^2,\;\;\;\forall v\in \rm{Sol}(D).
\end{equation}
The simplest example is of course
\begin{equation}
\label{def K0}
K[v]=\frac{1}{2}\int_{D} 
|v(x)|^2dx.
\end{equation}
Then, we are given a ``functional'' $\varphi\in E\rightarrow \mathcal{E}
[\varphi]\in\mathbf{R}$ on a suitable function space $E$,
the canonical example for us being the Dirichlet integral (\ref{dirichlet})
over the Sobolev space $E=H^1_0(D)$.
When we evolve $\varphi$ according to (\ref{transport}), we formally get
$$
\frac{d}{dt}\mathcal{E}[\varphi_t]=\int_D \mathcal{E}'[\varphi_t](x)\partial_t\varphi_t(x)dx
=-\int_D \mathcal{E}'[\varphi_t](x)\nabla\cdot(v_t(x)\varphi_t(x))dx,
$$
where we denote 
by $\mathcal{E}'$ the gradient of $\mathcal{E}$ with respect to the $L^2$ metric. (In the case
of the Dirichlet integral, $\mathcal{E}'[\varphi]=-\bigtriangleup\varphi$.)
Thus
\begin{equation}
\label{balance}
\frac{d}{dt}\mathcal{E}[\varphi_t]
=-\int_D \mathcal{E}'[\varphi_t](x)\nabla\varphi_t(x)
\cdot v_t(x)dx
\end{equation}
(using that $v_t$ is divergence-free).
\\
Then, we may write (\ref{balance}) as:
\begin{equation}
\frac{d}{dt}\mathcal{E}[\varphi_t]=-((G_t, v_t)),
\end{equation}
\begin{equation}
\label{G0}
G_t=\LERAY(\mathcal{E}'[\varphi_t]\nabla\varphi_t ),
\end{equation}
where we denote 
by $\LERAY$ 
the $L^2$ projection onto $\rm{Sol}(D)$.
Thus,
denoting by $K^*$ the Legendre-Fenchel transform $K^*[g]=\sup_w ((g,w))-K[w]$,
we get:
$$
\frac{d}{dt}\mathcal{E}[\varphi_t]+K[v_t]+K^*[G_t]=K[v_t]+K^*[G_t]-((G_t,v_t))
$$
where
by definition of the Legendre-Fenchel transform, 
the right-hand side 
is always nonnegative and vanishes
if and only if 
\begin{equation}
\label{closure}
v_t={K^*}'[G_t]
\end{equation}
(for instance, in case (\ref{def K0}), $v_t=G_t$).
Equation (\ref{closure}) precisely is the ``closure equation''
we need to define the "gradient flow" of $\mathcal{E}$ 
with respect to the evolution equation (\ref{transport}) with ``metric'' 
$K$.
(This way, we closely follow \cite{Bre}, in the spirit of \cite{AGS1,AGS2}.)
As just seen, this closure equation is equivalent to the differential inequality
$$
\frac{d}{dt}\mathcal{E}[\varphi_t]+K[v_t]+K^*[G_t]\le 0,
$$
or, using the definition of $K^*$ as the Legendre-Fenchel transform of $K$,
\begin{equation}
\label{variational0}
\frac{d}{dt}\mathcal{E}[\varphi_t]+K[v_t]+
((G_t,z_t))
-K[z_t]\le 0,
\;\forall\;\;z_t\in \rm{Sol}(D).
\end{equation}
So, our gradient flow is now defined by combining transport equation (\ref{transport}), definition 
(\ref{G0}), and either the closure equation (\ref{closure}) or the variational inequality (\ref{variational0}), which are
formally equivalent.

\section{The gradient flow equation}

From now on, let us consider, for simplicity, the case of
the periodic cube  $D=\mathbf{T}^d$, instead of a bounded domain of $\mathbf{R}^d$.
Accordingly,
all functions ($\varphi_t$, $v_t$, etc...)  to be considered will be of zero mean in $x\in D$.
We also concentrate on the case when:
\\
i) $\mathcal{E}$ is the Dirichlet integral (\ref{dirichlet});
\\
ii) $K$ is the Sobolev (semi-)norm of order $m$, $m\in\{0,1,2,\cdot\cdot\cdot\}$,
\begin{equation}
\label{def K}
K[v]=\frac{1}{2}\int_{D} 
|\nabla^m v(x)|^2dx,
\end{equation}
which should be understood, when $m>0$, as the $H^m$ Sobolev semi-norm of $v$ when it makes 
sense and $+\infty$ otherwise. 
\\
Then, $\mathcal{E}'=-\bigtriangleup$, ${K^*}'=(-\bigtriangleup)^{-m}.$
Also notice that the "relative entropy" reads, for each pair $(v,\omega)$ in $\rm{Sol}(D)$ with $\omega$ smooth,
$$
\eta_K[v,\omega]=K[v]-K[\omega]-((K'[\omega],v-\omega))=K[v-\omega]\ge c||v-\omega||^2,
$$
where $c>0$ depends only on $m$ and $d$, by Poincar\'e's inequality on the periodic cube.
So, equations (\ref{closure}) and (\ref{G0}) respectively become
$$
v_t=(-\bigtriangleup)^{-m}G_t\;,
$$
$$
G_t=\LERAY(\mathcal{E}'[\varphi_t]\nabla\varphi_t )=
\LERAY(-\bigtriangleup\varphi_t\nabla\varphi_t )=
\LERAY[-\nabla\cdot(\nabla \varphi_t \otimes \nabla \varphi_t)
+\frac{1}{2}\nabla(|\nabla\varphi_t|^2)]
$$
$$
=-\LERAY\nabla\cdot(\nabla \varphi_t \otimes \nabla \varphi_t)
$$
(since $\LERAY$ is the $L^2$ projection onto $\rm{Sol}(D)$ and, therefore, cancels any gradient).
Thus
\begin{equation}
\label{G}
v_t=-(-\bigtriangleup)^{-m}
\LERAY \nabla\cdot(\nabla \varphi_t \otimes \nabla \varphi_t).
\end{equation}
Similarly, (\ref{variational0}) becomes
\begin{equation}
\label{variational}
\frac{d}{dt}||\nabla\varphi_t||
^2+2K[v_t]
+((\nabla\varphi_t\otimes\nabla\varphi_t,\nabla z_t+\nabla z_t^T))\le 2K[z_t],
\end{equation}
for all smooth $z_t\in \rm{Sol}(D)$. 
\\
Finally, the gradient-flow equation reads:
\begin{equation}
\label{GF}
\partial_t \varphi_t+\nabla\cdot(\varphi_t v_t)=0,
\;\;\;(-\bigtriangleup)^m v_t=-\LERAY\nabla\cdot(\nabla\varphi_t\otimes\nabla\varphi_t).
\end{equation}
\subsection{
Physical interpretation of the GF equation
}
Physically speaking, 
the GF (gradient flow) equation (\ref{GF})
in the case $m=1$ corresponds to 
the ``Stokes flow''
\begin{equation}
\label{stokes}
\;\;\;\
\partial_t \varphi_t+\nabla\cdot(\varphi_t v_t)=0,
\;\;\;-\bigtriangleup v_t=-\LERAY\nabla\cdot(\nabla\varphi_t\otimes\nabla\varphi_t)
\end{equation}
of an electrically charged incompressible fluid (${v}$ and ${\varphi}$ being the velocity and the
electric potential), while the case $m=0$ rather corresponds to a ``Darcy flow''
\begin{equation}
\label{darcy}
\;\;\;\
\partial_t \varphi_t+\nabla\cdot(\varphi_t v_t)=0,
\;\;\;v_t=-\LERAY\nabla\cdot(\nabla\varphi_t\otimes\nabla\varphi_t).
\end{equation}
These are dissipative versions of the corresponding Euler equations
\begin{equation}
\label{euler}
\;\;\;\
\partial_t \varphi_t+\nabla\cdot(\varphi_t v_t)=0,
\;\;\;\partial_tv_t+\nabla\cdot(v_t\otimes v_t)=-\LERAY\nabla\cdot(\nabla\varphi_t\otimes\nabla\varphi_t).
\end{equation}

%%%%%%%%%%%%%%%%%%%%%%%%%%%%%%%%%%%%
\subsection{Special solutions and linear algebra
}
With a suitable potential added to the Dirichlet integral,
namely
\begin{equation}
\label{dirichlet pot}
\mathcal{E}[\varphi]=\frac{1}{2}\int_{D} (|\nabla\varphi(x)|^2-Qx\cdot x)dx,
\end{equation}
where $Q$ is a fixed $d\times d$ symmetric matrix,
and set on the unit ball instead of a periodic box,
the gradient-flow (GF) equation 
has interesting special solutions which are linear in ${x}$
\begin{equation}
\label{}
\;\;\;\
\nabla\varphi_t(x)=M_tx,\;\;\;v_t(x)=V_tx,\;\;\;M_t=M_t^T,\;\;\;V_t=-V_t^T.
\end{equation}
The resulting equation reads
\begin{equation}
\label{brockett}
\;\;\;\
\frac{dM_t}{dt}=[V_t,M_t],\;\;\;V_t=[M_t,Q].
\end{equation}
With (\ref{brockett}), we recover
the Brockett diagonalizing gradient flow for $d\times d$ symmetric matrices 
(recently revisited by Bach and Bru, in the generalized case of
infinite dimensional for self-adjoint operators) \cite{Bro,BaBr}
(see also \cite{GaHo} in connection with Fluid Mechanics).
The case when $Q={\rm{diag}}(1,2,\cdot\cdot\cdot,d)$
is of peculiar interest. In that case, $M_t$ converges to its diagonal form
(with eigenvalues sorted in non-decreasing order) as $t$ goes to $+\infty$.

\section{Analysis of the gradient flow equation}

The last part of this article is devoted to the analysis of the gradient flow equation \ref{GF}.
For that purpose, we closely follow the ideas and concepts of our recent work \cite{Bre}.

\subsection{A concept of ``dissipative solutions''}

From the analysis viewpoint, 
we ignore whether or not gradient-flow (GF) equation (\ref{GF}), namely 
$$
\partial_t \varphi_t+\nabla\cdot(\varphi_t v_t)=0,
\;\;\;(-\bigtriangleup)^m v_t=-\LERAY\nabla\cdot(\nabla\varphi_t\otimes\nabla\varphi_t)\;,
$$
is locally well-posed
in any space of smooth functions (unless $m>d/2+1$). 
The global existence of weak solutions can be expected for
the Stokes version (\ref{stokes}) (with $\varphi$ a priori
in $L^\infty_t(H^1_x)$
and $v$ ``almost'' in $L^\infty_t(W^{1,1}_x)$), while such a result 
looks out of reach in the case of the ``Darcy'' version (\ref{darcy}).
\\
Anyway, we prefer a much more
"robust" concept of solutions, that we call ``dissipative''
after \cite{Bre},
somewhat in the spirit of Lions' dissipative solutions to the Euler equations \cite{Li} 
and following some ideas of the analysis of the linear heat equations for general measured
metric spaces by Ambrosio, Gigli and Savar\'e \cite{AGS2}.
We keep transport equation (\ref{transport}) and integrate 
(\ref{variational}) on $[0,t]$, for all $t\ge 0$, with a suitable
exponential weight, which leads to:
\begin{equation}
\label{dissipative}
\int_0^t \{
2K[v_s]+
((\nabla\varphi_s\otimes\nabla\varphi_s,rI_d+\nabla z_s+\nabla z_s^T))
-2K[z_s]\}e^{-s r}ds
\end{equation}
$$
+||\nabla\varphi_t||^2e^{-tr}  \le ||\nabla\varphi_0||^2,
 \;\;\;
 {\rm{for\;every\;smooth\;field\;}}
 t\rightarrow z_t\in \rm{Sol}(D).
$$
Here $r\ge 0$ is a constant, depending on $z$, chosen so that
\begin{equation}
\label{positive}
\forall (t,x),\;\;rI_d+\nabla z_t(x)+\nabla z_t(x)^T\ge 0,\;{\rm{\;in\;the\;sense\;of\;symmetric\;matrices,}}
\end{equation}
in order to be sure that inequality (\ref{dissipative}) only involves convex functionals of $(\varphi,v)$.
From the functional analysis viewpoint, it is natural to consider solutions 
$(B_t=\nabla\varphi_t,v_t,t\in [0,T])$, for each fixed $T>0$, in the space
$$
C^0_w([0,T],L^2(D,\mathbf{R}^2))\times L^2([0,T],{\rm{Sol}}(D)),
$$
where $C^0_w(L^2)$ just means continuity in time with respect to the weak topology of $L^2$.
\newpage
%%%%%%%%%%%%%%%%%%%%%%%%%%%%%%%%%%%%
\subsection{Uniqueness of smooth solutions among dissipative solutions}
%%%%%%%%%%%%%%%%%%%%%%%%%%%%%%%%%%%%

\begin{theorem}
\label{uniqueness}
Assume $D=(\ERRE/\ZED)^d$ and define $K$ by (\ref{def K}), namely
$$
K[v]=\frac{1}{2}\int_{D} 
|\nabla^m v(x)|^2dx,
$$
with "relative entropy"
\begin{equation}
\label{entropy1}
\eta_K[a,b]=K[a]-K[b]-((K'[b],a-b))=
\frac{1}{2}\int_{D} 
|\nabla^m (a-b)(x)|^2dx.
\end{equation}
Let us fix $T>0$ and consider
$$
(B_t=\nabla\varphi_t,v_t,t\in [0,T])\in C^0_w([0,T],L^2(D,\mathbf{R}^2))\times L^2([0,T],{\rm{Sol}}(D)),
$$
a dissipative solution of the GF equation (\ref{GF}) up to time $T$,
in the sense of (\ref{transport},\ref{dissipative},\ref{positive}). Let 
$(\beta_t=\nabla\psi_t,\omega_t,t\in [0,T])$ 
be any pair of smooth functions with $\omega$ valued 
in $\rm{Sol}(D)$. Then there is a constant $C$ depending only on $K$ and
the spatial Lipschitz constant of $(\beta,\omega)$, up to time $T$,
so that, for all $t\in [0,T]$,
\begin{equation}
\begin{split}
||B_t-\beta_t||^2
+\int_0^t e^{(t-s)C}\{\eta_K[v_s,\omega_s]ds-2J_s^{L}\}ds
\le ||B_0-\beta_0||^2e^{tC}
\\
J_t^L=-((B_t-\beta_t,\nabla(\omega_t\cdot\beta)+\partial_t \beta_t))
-((\LERAY\nabla\cdot(\beta_t\otimes \beta_t)+K'[\omega_t],v_t-\omega_t)).
\end{split}
\end{equation}
In particular, $J_t^L$ exactly vanishes as
$(\beta=\nabla\psi,\omega)$ is a smooth solution to the GF equation (\ref{GF}), namely
$$
\partial_t \beta_{t}+\nabla(\omega_{t}\cdot\beta_{t})=0,
\;\;\;K'[\omega_t]=-\LERAY\nabla(\beta_t\otimes\beta_t),
$$
in which case
\begin{equation}
\label{strong uniqueness}
\begin{split}
||B_t-\beta_t||^2
+ \int_0^t e^{-(t-s)C}\eta_K[v_s,\omega_s]ds\;\le 
||B_0-\beta_0||^2e^{-tC}.
\end{split}
\end{equation}
\end{theorem}
This implies the uniqueness of smooth
solutions among all dissipative solutions, for any given prescribed 
smooth initial condition.

%%%%%%%%%%%%%%%%%%%%%%%%%%%%%%%%%%%%%%%%%
\subsection{Global existence of dissipative solutions}
%%%%%%%%%%%%%%%%%%%%%%%%%%%%%%%%%%%%%%%%%

In the spirit of \cite{Bre}, 
at least in the case:  $D=\mathbf{T}^d$,
$$
\mathcal{E}[\varphi]=\frac{1}{2}\int_{\mathbf{T}^d} |\nabla\varphi(x)|^2dx,\;\;\;
K[v]=\frac{1}{2}\int_{\mathbf{T}^d} 
|(\nabla)^m v(x)|^2dx\;\;\;(m=0,1,2,\cdot\cdot\cdot),
$$
it is fairly easy to establish, for the "dissipative"  formulation (\ref{transport},\ref{dissipative},\ref{positive}) of (\ref{transport},\ref{closure})
and
for each initial condition $\varphi_0$ with finite Dirichlet integral,
the existence of a global solution $(B=\nabla\varphi,v)$ in
$
C_w^0(\mathbf{R}_+,L^2(D,\mathbf{R}^d))\times L^2(\mathbf{R}_+,{\rm{Sol}}(D)).
$
Without entering into details, let us sketch the proof.
We approximate $B_0$ strongly in $L^2$ by some smooth field 
$B_0^\epsilon=\nabla\varphi_0^\epsilon$
and mollify $K$ by substituting for it
$$
K^{M,\epsilon}(v)=K(v)+\epsilon ||\nabla^M v||^2
$$
with $M$ sufficiently large ($M>d/2+1$) and $\epsilon>0$. In this case, we get, $M$ and $\epsilon$ being fixed,
a uniform a priori bound for $v$ 
in $L^2([0,T],C^1(D))$, which is enough to solve transport equation (\ref{transport}) in the classical framework of the
Cauchy-Lipschitz theory of ODEs. Then, we get a smooth approximate solution $(B^\epsilon=\nabla \varphi^\epsilon,v^\epsilon)$
satisfying transport equation (\ref{transport}), i.e.,
$$
\partial_t\varphi_t^\epsilon
+\nabla\cdot(\varphi_t^\epsilon v_t^\epsilon)=0,
$$
together with (\ref{dissipative},\ref{positive}), namely
$$
\int_0^t \{
2K^{M,\epsilon}
[v^\epsilon_s]+
((B^\epsilon_s\otimes B^\epsilon_s,rI_d+\nabla z_s+\nabla z_s^T))
-2K^{M,\epsilon}[z_s]\}e^{-s r}ds
$$
$$
+||B^\epsilon_t||^2e^{-tr}  \le ||B^\epsilon_0||^2,
 \;\;\;
 {\rm{for\;every\;smooth\;field\;}}
 t\rightarrow z_t\in \rm{Sol}(D)
$$
satisfying (\ref{positive}), and, in particular (for $z=0$)
$$
\int_0^t 2K^{M,\epsilon}
[v^\epsilon_s]ds
+||B^\epsilon_t||^2 \le ||B^\epsilon_0||^2.
$$
We get enough compactness for the approximate solutions to get a limit $(B,v)$
in space 
$$
C_w^0(\mathbf{R}_+,L^2(D,\mathbf{R^d}))\times L^2(\mathbf{R}_+,{\rm{Sol}}(D)),
$$
and pass to the limit in the transport equation (since $B^\epsilon=\nabla \varphi^\epsilon$).
Finally, by lower semi-continuity, we may pass to the limit in the dissipation inequality and
obtain (\ref{dissipative},\ref{positive}), which concludes the (sketch of) proof.
\\
Observe that, for $m\ge 1$, the $L^2$ norm of $\nabla v_t$ is square integrable in time. This implies
by DiPerna-Lions ODE theory (see \cite{Li}), as already discussed, 
that the law of $\varphi_t$
stays unchanged during the evolution by (\ref{transport})
(but, unless $m>1+d/2$, not necessarily its topology, which is of some interest in view of the minimization problem
(\ref{kelvin}) we started with).
However, unless $m>1+d/2$, it is unclear to us that (\ref{transport},\ref{closure}) even admits local smooth solutions.

%%%%%%%%%%%%%%%%%%%%%%%%%%%%%%%%%%%%%%%%%
%%%%%%%%%%%%%%%%%%%%%%%%%%%%%%%%%%%%%%%%%
\section*{Appendix: Proof of Theorem \ref{uniqueness}}
Choose $r\ge 0$ such that $\omega$ satisfies (\ref{positive}), namely
$$
\forall (t,x),\;\;rI_d+ \nabla \omega_t(x)+\nabla \omega_t(x)^T\ge 0,\;{\rm{\;in\;the\;sense\;of\;symmetric\;matrices.}}
$$
Since $(B=\nabla\varphi,v)$ is a dissipative solution, we get, by setting $z=\omega$ in definition (\ref{dissipative}),
\begin{equation}
\int_0^t \{
2K[v_s]+
((B_s\otimes B_s,rI_d+\nabla \omega_s+\nabla \omega_s^T))
-2K[\omega_s]\}e^{-s r}ds
\end{equation}
$$
+||B_t||^2e^{-tr}  \le ||B_0||^2.
$$
Let us now introduce for each $t\in [0,T]$
\begin{equation}
\label{Nt}
N_t=||B_0||^2 e^{rt}-\int_0^t 
\{2K[v_s]+
((B_s\otimes B_s,rI_d+\nabla \omega_s+\nabla \omega_s^T))
-2K[\omega_s]\}e^{r(t-s)}ds
\end{equation}
so that
$$
N_t\ge ||B_t||^2, \forall t\in [0,T].
$$
By definition (\ref{Nt}) of $N_t$, we have
$$
(\frac{d}{dt}-r)N_t=-2K[v_t]-((B_t\otimes B_t,rI_d+\nabla \omega_t+\nabla \omega_t^T))+2K[\omega_t]
$$
and, therefore,
\begin{equation}
\label{dtNt}
\frac{d}{dt}N_t
=r(N_t-||B_t||^2)-2K[v_t]+((B_t\otimes B_t,\nabla \omega_t+\nabla \omega_t^T))+2K[\omega_t]
\end{equation}
(in the distributional sense and also for a.e. $t\in [0,T]$).
\\
We now want to estimate
\begin{equation}
\label{et}
e_t=
||B_t-\beta_t||^2+(N_t-||B_t||^2)
=N_t-2((B_t,\beta_t))+||\beta_t||^2,\;\forall t\in [0,T],
\end{equation}
where $\beta_t=\nabla\psi_t$.
Since $(B=\nabla\varphi,v)$ is a dissipative solution, it solves transport equation (\ref{transport})
which implies
$$
\partial_t B_t+\nabla(B_t\cdot v_t)=0,
$$
after derivation in $x$.
Thus
$$
\frac{d}{dt}((B_t,\beta_t))
=\int B_{ti}(\beta_{ti,t}+v_{ti}\beta_{tj,j})
$$
(where we use notations $\beta_{ti,j}=\partial_j(\beta_t)_{i}$, etc...and skip
summations on repeated indices $i$, $j$...).
\\
Using (\ref{dtNt}) and definition (\ref{et}), we deduce
$$
\frac{d}{dt}e_t
=r(N_t-||B_t||^2)-2K[v_t]-((B_t\otimes B_t,\nabla \omega_t+\nabla \omega_t^T))+2K[\omega_t]
$$
$$
+\int 2(\beta-B)_{ti}\beta_{ti,t}-2B_{ti}v_{ti}\beta_{tj,j}.
$$
Thus
$$
\frac{d}{dt}e_t=
r(N_t-||B_t||^2)-2K[v_t]+2K[\omega_t]+J_t
$$
where
$$
J_t=\int -B_{ti}B_{tj}(\omega_{ti,j}+\omega_{tj,i})
+2(\beta-B)_{ti}\beta_{ti,t}-2B_{ti}v_{ti}\beta_{tj,j}.
$$
Denoting the "relative entropy" of $K$ by  $\eta_K[a,b]=K[a]-K[b]-((K'[b],a-b))$, we have obtained
\begin{equation}
\label{dtet}
\frac{d}{dt}e_t+2\eta_K[v_t,\omega_t]=
r(N_t-||B_t||^2)+J_t-2((K'[\omega_t],v_t-\omega_t))
\end{equation}
We may write
$$
J_t=
J_t^{Q}+J_t^{L1}+J_t^{L2}+J_t^{C}
$$
where 
$J_t^{Q}$, $J_t^{L1}$, $J_t^{L2}$, $J_t^{C}$ are respectively 
quadratic, linear, linear, and constant with respect to $B-\beta$ and $v-\omega$,
with coefficient depending only on $\omega,\beta$:
$$
J_t^{Q}=\int -(B-\beta)_{ti}(B-\beta)_{tj}(\omega_{ti,j}+\omega_{tj,i})
-2(B-\beta)_{ti}(v-\omega)_{ti}\beta_{tj,j}
$$
$$
J_t^{L1}=\int 2(B-\beta)_{ti}[-\beta_{tj}(\omega_{ti,j}+\omega_{tj,i})-\beta_{ti,t}
-\omega_{ti}\beta_{tj,j}]
$$
$$
J_t^{L2}=-\int 2(v-\omega)_{ti}\beta_{ti}\beta_{tj,j}
$$
$$
J_t^{C}=\int [-\beta_{ti}\beta_{tj}(\omega_{ti,j}+\omega_{tj,i})
-2\beta_{ti}\omega_{ti}\beta_{tj,j}].
$$
%%%%%%%%%%%%%%%%%%%%%%%%%%%%%%%%%%%%%%%
Let us reorganize these four terms.
By integration by part of its first term, we see that $J_t^{C}=0$, using
that $\beta$ is a gradient and $\omega$ is divergence-free. 
More precisely
$$
-\int \beta_{ti}\beta_{tj}(\omega_{ti,j}+\omega_{tj,i})=-2\int \beta_{ti}\beta_{tj}\omega_{ti,j}
=2\int \beta_{ti}\beta_{tj,j}\omega_{ti}+2\int \beta_{ti}\beta_{tj,j}\omega_{ti}
$$
$$
=0+2\int \beta_{ti}\beta_{tj,j}\omega_{ti}.
$$
Using that $v_t-\omega_t$ is divergence-free while $\beta$ is a gradient, we immediately get
$$
J_t^{L2}=-2((\LERAY\nabla\cdot(\beta_t\otimes \beta_t),v_t-\omega_t)).
$$
Next, 
we find
$$
J_t^{L1}=-2((B_t-\beta_t,\nabla(\omega_t\cdot\beta)+\partial_t \beta_t)).
$$
Indeed, since
$$
J_t^{L1}=\int 2(B-\beta)_{ti}[-\beta_{tj}(\omega_{ti,j}+\omega_{tj,i})-\beta_{ti,t}
-\omega_{ti}\beta_{tj,j}]
$$
we get
$$
J_t^{L1}+2((B_t-\beta_t,\nabla(\omega_t\cdot\beta)+\partial_t \beta_t))
=\int 2(B-\beta)_{ti}[-\beta_{tj}(\omega_{ti,j}+\omega_{tj,i})+(\omega_{tj}\beta_{tj}),_i
-\omega_{ti}\beta_{tj,j}]
$$
$$
=\int 2(B-\beta)_{ti}[-(\beta_{tj}\omega_{ti})_j-\beta_{tj}\omega_{tj,i}+\omega_{tj,i}\beta_{tj}+\omega_{tj}\beta_{tj,i}]
=\int 2(B-\beta)_{ti}[-(\beta_{tj}\omega_{ti})_j+\omega_{tj}\beta_{tj,i}]
$$
$$
=\int 2(B-\beta)_{ti}[-(\beta_{tj}\omega_{ti}),_j+(\beta_{ti}\omega_{tj}),_j-\beta_{ti}\omega_{tj,j}]=0
$$
(since $\omega$ is divergence-free while $B_t-\beta_t$ and $\beta_t$ are gradients).
Next, since
$$
J_t^{Q}=\int -(B-\beta)_{ti}(B-\beta)_{tj}(\omega_{ti,j}+\omega_{tj,i})
-2(B-\beta)_{ti}(v-\omega)_{ti}\beta_{tj,j}
$$
we may find,
for any fixed $\epsilon>0$, 
a constant $C_\epsilon$ 
(depending on the spatial Lipschitz constant of $(\beta,\omega)$)
such that
$$
J_t^{Q}\le \epsilon ||v_t-\omega_t||^2+C_\epsilon ||B_t-\beta_t||^2.
$$
Using (\ref{entropy1}), we may choose $\epsilon$ small enough so that
$$
\epsilon ||v_t-\omega_t||^2\le \eta_K[v_t,\omega_t].
$$
So, we get from (\ref{dtet})
\begin{equation}
\label{dtet2}
\frac{d}{dt}e_t+\eta_K[v_t,\omega_t]\le r(N_t-||B_t||^2)
+C_\epsilon ||B_t-\beta_t||^2
+2J_t^L
\end{equation}
where
$$
J_t^L=-((B_t-\beta_t,\nabla(\omega_t\cdot\beta)+\partial_t \beta_t))-((\LERAY\nabla\cdot(\beta_t\otimes \beta_t)
+K'[\omega_t],v_t-\omega_t)).
$$
By definition (\ref{et}) of $e_t$, namely
$
e_t=
||B_t-\beta_t||^2+(N_t-||B_t||^2),
$
we have obtained 
$$
\frac{d}{dt}e_t+\eta_K[v_t,\omega_t]\le Ce_t 
+2J_t^L
$$
for a constant $C$ depending only on $\beta$, $\omega$ and $K$.
By integration we deduce
$$
e_t+\int_0^t e^{(t-s)C}\eta_K[v_s,\omega_s]ds
\le e_0e^{tC}+2\int_0^t e^{(t-s)C}J_s^{L}ds.
$$
Next, let us remind that
$e_t\ge ||B_t-\beta_t||^2$ with equality for $t=0$
(since $N_t\ge ||B_t||^2$ with equality at $t=0$).
Thus, we have finally obtained
$$
||B_t-\beta_t||^2
+\int_0^t e^{(t-s)C}\eta_K[v_s,\omega_s]ds
\le ||B_0-\beta_0||^2e^{tC}+2\int_0^t e^{(t-s)C}J_s^{L}ds
$$
with
$$
J_t^L=-((B_t-\beta_t,\nabla(\omega_t\cdot\beta)+\partial_t \beta_t))-((\LERAY\nabla\cdot(\beta_t\otimes \beta_t)
+K'[\omega_t],v_t-\omega_t))
$$
and the proof of Theorem \ref{uniqueness} is now complete.

\subsection*{Acknowlegments}
This  work has been partly supported by the ANR contract ISOTACE
(ANR-12-MONU-013).
\\
The author thanks M. Salmhofer for pointing out references
\cite{Bro,BaBr}.

%%%%%%%%%%%%%%%%%%%%%%%%%%%%%%%%%%%%%%%%%%%%%%%%%%%%%%%%%%%%

\end{document}